\documentclass[12pt]{article}
\usepackage{amssymb}

\newcommand\eps{\varepsilon}
\newcommand\ibc{\textsc{ibc}}
\newcommand\ivp{\textsc{ivp}}
\newcommand\bvp{\textsc{bvp}}
\newcommand\ode{\textsc{ode}}

\begin{document}

\title{Open letter on ``Adaptivity and computational complexity 
in the numerical solution of ODEs'' 
by
Silvana Ilie, 
Gustaf S\"oderlind and 
Robert M. Corless}

\author{Erich Novak \\
Mathematisches Institut, Universit\"at Jena\\
Ernst-Abbe-Platz 2, 07740 Jena, Germany\\
email:\ novak@mathematik.uni-jena.de\\
\qquad
\\
Henryk Wo\'zniakowski \\ 
Department of Computer Science, Columbia University,\\
New York, NY 10027, USA, and\\
Institute of Applied Mathematics, University of Warsaw\\
ul. Banacha 2, 02-097 Warszawa, Poland\\
email:\ henryk@cs.columbia.edu}

\maketitle

\begin{abstract} 
This is an open 
letter on the paper  ``Adaptivity and computational complexity 
in the numerical solution of ODEs'' 
by
Silvana Ilie, 
Gustaf S\"oderlind and 
Robert M. Corless. 
We sent this letter to the authors in August 2008. 
\end{abstract} 

\noindent
Dear Silvana, Gustaf and Robert,

Your paper ``Adaptivity and computational complexity in the numerical
solution of \ode{}s'' recently appeared in the \emph{Journal of
Complexity}.  After 
studying this paper, we found a few points
that we would like to discuss with you.

In your paper you study the numerical solution of \ivp{}s and \bvp{}s
for \ode{}s.  In the abstract you claim that (under some mild
conditions) ``an adaptive method is always more efficient than a
non-adaptive method''.

You mention ``computational complexity'' in the title and one of your
keywords is ``information-based complexity''.  Moreover, you cite the
book of Werschulz~\cite{We91}.  It seems therefore that your paper is
meant as a contribution to complexity theory and, more particular, to
\ibc.

The first two sentences of the introduction read as follows: 

\begin{quote}
  The complexity of numerical algorithms is central to the assessment
  of computational performance. For some algorithms, like in linear
  algebra, the complexity is well known and established; for others,
  like ordinary differential equations (\ode{}s), the complexity is still
  open to analysis.
\end{quote} 

We agree with the first sentence.  However, you are conflating two
different concepts: the \emph{cost} of an algorithm and the
\emph{complexity} of a problem.  The complexity of a problem is the
cost of an \emph{optimal} algorithm.  In other words, the complexity
is the minimal computational cost of solving the problem, said minimum
being taken over \emph{all} algorithms that solve the problem.
\emph{Complexity is \textbf{not} a property of any particular algorithm;
it is a property of the problem under consideration.}

More important is the following.  We do not currently know the
complexity of many problems of linear algebra. Probably the most basic
problem from linear algebra is 
 matrix multiplication, which is
equivalent to solving linear equations. Unfortunately, we do
\emph{not} know the complexity of this problem. We know only bounds on
the complexity. This  famous problem has its own history, which has
been enriched by the contributions of Strassen, Pan, Coppersmith and
Winograd. They would be quite surprised to see that someone claims
that the complexity of this problem is well known.

On the other hand, if we interpret complexity as being equivalent to
cost, then we agree that the cost of some algorithms in linear algebra
is known.  But we also know the cost of many other algorithms for
\ode{}s, such as Euler, Runge-Kutta, and multistep methods.  This
suggests that the notion of complexity used by you is different than
the notion of cost, as it should be.

Of course we believe that the complexity of \ode{}s ``is still open to
analysis'' since we do not know \emph{everything} concerning the
complexity of \ode{}s.  Nevertheless, it is quite surprising for us
that you almost completely ignore what we \emph{do} know about the
complexity of solving \ode{}s.  Kacewicz proved many complexity
results for \ode{}s, mainly for \ivp{}s. We recommend his
paper~\cite{Ka87} on the complexity of \ivp{}s for \ode{}s and stress
that most of the upper bounds are proved using \emph{adaptive}
algorithms.  Moreover, Kacewicz wrote another paper~\cite{Ka90}, in
 which he compares the power of adaptive and non-adaptive information
and proves that adaption is \emph{exponentially} better for
$d$-variate \ode{}s.

\emph{There is no \ibc{} result saying that adaption \textbf{never}
helps}, see also \cite{No96,TWW88}.  Werschulz proved many
complexity results for \ode{}s, mainly for \bvp{}s.  Other authors
proved many more results.  None of these results can be found in your
article, although the book of Werschulz~\cite{We91} is cited.

You probably heard the result that adaption does \emph{not} help but
this result, as all mathematical results, only holds under specific
hypotheses.  In this case, we must assume that 
\begin{itemize}
\item we want to approximate a \emph{linear} operator 
\item defined on a \emph{convex}, \emph{balanced} class, and that
\item we use linear functionals for approximating the linear operator.
\end{itemize}
Then adaption can be only better by a factor of at most $2$.  Note
that the operator associated with \ode{}s is \emph{non-linear},
since the solution of \ode{}'s of the form $u'(t)=f(u(t),t)$ depends
non-linearly on~$f$.  So, the theorem on adaption cannot be applied
for \ode{}s.  In fact, the results in Kacewicz~\cite{Ka90} 
dramatic prove that adaption can be significantly better than
non-adaption. 

Your introduction continues:

\begin{quote}
  In the former case, the problems are ``computable'', meaning that
  (theoretically) the exact solution can be obtained after a finite
  number of operations, and this operation count then becomes a
  measure of the complexity. By contrast, for problems in analysis we
  can only compute approximate solutions converging to the exact
  solution. This makes an assessment of complexity more difficult, as
  the algorithmic complexity will depend on problem characteristics as
  well as the requested accuracy.
\end{quote} 

Of course, the necessity of computing approximations, rather than
exact solutions, is not restricted to problems (such as integration)
from analysis.  There are also problems from linear algebra (such as
eigenproblems) where we only can compute an approximation of the
solution, despite the fact that we have complete information.  Hence
it doesn't matter whether you're talking about analysis or linear
algebra; there are times when we need to compute an approximation
whose error (in some norm) is bounded by some positive~$\eps$.
This~$\eps$ depends on the application, and might be very small or
relatively large.

Your introduction continues further:

\begin{quote}
  In differential equations, adaptive algorithms are of fundamental
  importance.  Such algorithms attempt to minimize some (usually
  coarse) measure of complexity, subject to a prescribed accuracy
  criterion and the problem properties encountered during the
  computation. This is generally done by using \emph{non-uniform
  discretization grids} in order to put the discretization points
  where they matter most to accuracy, while keeping their total number
  small.

  Naturally, for some problems, uniform grids might be optimal from
  the point of view of complexity, e.g., if one considers FFT-based
  algorithms for Poisson's equation on a rectangular domain. For
  linear problems, similar considerations led Werschulz to question
  whether adaptive methods are more efficient, using a topological
  argument to show that the efficiency gain would be limited to a
  factor of two \cite[pp.~38--39]{We91}. In this paper, however, we
  will prove that adaptivity is better than non-adaptivity.
\end{quote} 

You never explain what you mean by ``some measure of complexity''.
Moreover, you also do not define your notion of adaptivity.  One has
to read between the lines to figure out what this really means; it
seems to us that you use the word \emph{adaption} to characterize
methods that are based on a \emph{non-uniform grid}.

First of all, we admit that any of us are free to use the term
``adaption'' to mean whatever we want. However, if someone wants to
compare her/his results with the \ibc{} literature, it seems necessary
to use the same definition of adaption as is used in \ibc, or at least
to point out differences.  Since your paper suggests that the 
\ibc{} result about adaption is wrong, it is a pity that your paper
does not precisely formulate this result.  Note that this result has a
venerable history.  Back in 1971, Bakhvalov~\cite{Ba71} proved a
version of this result for the approximation of linear functionals.
The result was proved in full generality for linear operators
independently by Gal and Micchelli~\cite{GM80} and Traub and
Wo\'zniakowski~\cite{TW80}.  This last result is the one mentioned by
Werschulz~\cite[pp.~38--39]{We91}.

Since your paper claims that the \ibc{} result on adaption is wrong,
you should indicate an error in the proof or a counterexample.  At
the very minimum, you should state this wrong result precisely.  Not
surprisingly, we think that the \ibc{} result on adaption is perfectly
fine.  It would be very surprising if someone would find an error in
such an old and elementary result.  After all
it has been checked by many people.  We are also
surprised that even today, after 27 years, the precise statement of
this result on adaption is not better known outside the \ibc{} community.

\emph{The notion of adaption used in \ibc{} has nothing to do with
uniform or non-uniform grids.}  Information is non-adaptive if we
use, say, function values at sample points that are the \emph{same}
for all functions from a given class. They can be from uniform or
non-uniform grid or from any set. Information is adaptive if sample
points vary with functions from a given set. So Gauss quadrature uses
non-adaptive information for univariate integration for a class of
smooth functions, and Newton, secant or bisection uses adaptive
information for solving univariate non-linear equations for
appropriate classes of functions.

We want to indicate
that your paper does \emph{not} contribute to the complexity
theory for \ode{}s. In the introduction, you say:

\begin{quote} 
  In this paper we analyze the complexity of solving \ode{}s using
  adaptive one-step methods based on local error control.
\end{quote} 

Here you admit that you only want to study a very specific class of
algorithms. In complexity theory, however, we want to discuss
\emph{all} algorithms and, in particular, find \emph{optimal}
algorithms.  If we restrict  the class of algorithms to ``one-step
methods based on local error control'',  then we cannot expect to
obtain sharp complexity results.

In Section~2 you present the problem as follows: 

\begin{quote}
  We shall consider the problem of solving an \ode{}, written as an
  operator equation $\mathfrak L(u)=f$ with either initial or
  boundary conditions.
\end{quote} 

You do \emph{not} introduce a \emph{class} of such problems.  Such a
class might be defined by specifying certain properties (say,
smoothness) of the functions~$f$ for which we wish to solve the
problem $\mathfrak L(u)=f$. The complexity depends on the domain of a
particular problem, and it may vastly change if we change the
domain. You only discuss a single \ode{}. There is no such a thing as
a complexity theory for the computation of a single object.  For any
kind of complexity theory one has to describe a \emph{class of
  problems} under consideration.  In the paper you prove results, in
particular about ``optimal'' grids.  A typical result, such as
Theorem~3, starts with
\begin{quote} 
  Let $p$ be the order of the method, $\dots$ 
\end{quote} 
You then \emph{assume} that a certain method is used and only 
discuss variants of the same method with different grid points. 
Again, such statements are a contribution to the cost analysis of a
particular method, but
not to complexity theory. After all, you do not even aim 
to present optimality results for the class of \emph{all algorithms}. 
You do \emph{not} prove a result that your  method is better than a
suitable method using \emph{non-adaptive} grid which would be
needed for a complexity result. 

Theorem 5 is about the ``minimum number of grid points'' and hence
looks like a complexity result. Again, it is not since you restrict
yourselves to the study of a very specific class of algorithms, which
is applied to a single operator equation. In a way, you admit this
shortcoming on p.~354 when you say that Theorem~5 is about the cost to
solve the \ode{} ``with the given method''.

Complexity theory is different.  We do not want to solve a single
problem instance ``with the given method''. The aim of complexity
theory is to construct and to analyze \emph{optimal} algorithms for
the solution of \emph{classes} of problems.

On p.~355, you claim:

\begin{quote} 
  This supports the ``conventional wisdom'' and resolves the complexity 
  controversy, \cite[p.~124]{We91}. 
\end{quote} 

Unfortunately, the adaption problem is not even touched in your paper.
Obviously the question ``For which problems are adaptive algorithms
superior?'' is very important.  Some answers to this question, as well
as exact definitions, can be found in the texts that we cite in this
letter.


We hope to get your comments on the points we raised in this letter.

\medskip

\begin{flushright}
  \begin{minipage}{0.6\linewidth}
    With best regards,

    \bigskip

    Erich Novak 

    Henryk Wo\'zniakowski


  \end{minipage}
\end{flushright}

\end{document}